\documentclass[aps,pre,reprint,groupedaddress,amsmath,amssymb]{revtex4-2}
\usepackage{graphicx}
\usepackage{bm}
\usepackage{color}

\begin{document}

\title{Connection among stochastic Hamilton-Jacobi-Bellman equation, path-integral, and Koopman operator on nonlinear stochastic optimal control}

\author{Jun Ohkubo}
\email[]{johkubo@mail.saitama-u.ac.jp}
\affiliation{
Graduate School of Science and Enginnering, Saitama University,
255 Shimo-Okubo, Sakura-ku, Saitama, 338-8570, Japan
}
\affiliation{
JST, PREST, 4-1-8 Honcho, Kawaguchi, Saitama, 332-0012, Japan
}


\begin{abstract}
The path-integral control, which stems from the stochastic Hamilton-Jacobi-Bellman equation, is one of the methods to control stochastic nonlinear systems. This paper gives a new insight into nonlinear stochastic optimal control problems from the perspective of Koopman operators. When a finite-dimensional dynamical system is nonlinear, the corresponding Koopman operator is linear. Although the Koopman operator is infinite-dimensional, adequate approximation makes it tractable and useful in some discussions and applications. Employing the Koopman operator perspective, it is clarified that only a specific type of observable is enough to be focused on in the control problem. This fact becomes easier to understand via path-integral control. Furthermore, the focus on the specific observable leads to a natural power-series expansion; coupled ordinary differential equations for discrete-state space systems are derived. A demonstration for nonlinear stochastic optimal control shows that the derived equations work well.
\end{abstract}


\maketitle

\section{Introduction}
\label{sec_introduction}

There are many stochastic phenomena in the real world, and the studies for them are interesting topics in physics. Furthermore, the focus includes not only analyzing the stochastic phenomena but also controlling them. One of the famous examples is the stochastic path-integral control by Kappen \cite{Kappen2005,Kappen2005a}. In control problems for deterministic cases, the Hamilton-Jacobi-Bellman (HJB) equation is usually employed \cite{Stengel_book}. The HJB equation is also useful for the cases with stochasticity, and in Kappen's prominent papers \cite{Kappen2005,Kappen2005a}, the usage of the path-integral approach and importance sampling methods gives a novel way for the control problem from the viewpoint of physics. There are many papers on path-integral control, and one of the improvements is the iterative path-integral method derived from information-theoretic arguments \cite{Williams2016,Williams2017,Williams2017a,Williams2017b}. In the machine learning community, a connection to the gradient descent method makes more practical algorithms \cite{Okada2018}. As shown in these papers, the iterative path-integral methods give real-time control, and there are several practical applications. Although the usage of neural networks is one of the promising approaches for control issues, it gives black-boxes in general. In physics, there are many theoretical frameworks for dealing with stochastic systems, and it is useful to explore new methods for control problems from a physical perspective. Such studies will avoid black-box methods for control issues.

One of the candidates to seek such control methods is the Koopman operator \cite{Mauroy2020}. The Koopman operator is a kind of composition operator and typically applied to dynamical systems. Although the original dynamical systems are nonlinear and finite-dimensional, the Koopman operator gives a map from functions to functions; the operator is linear but infinite-dimensional, as denoted later. Because the linearity is tractable in many situations, the Koopman operator has recently attracted much attention in time-series data analysis \cite{Mezic2005,Williams2015}. In such time-series data analysis, the dynamic mode decomposition gives much insight for various topics ranging from fluid dynamics to nonlinear oscillating systems. Because of its compatibility with dynamical systems, the Koopman operator is also used in control problems \cite{Peitz2019,Li2019,Klus2020,Korda2020,Peitz2020}. However, there is much room for research on the applications of the Koopman operator in control issues.

In many cases, a kind of sampling method gives the evaluation of statistics. However, in the path-integral or the iterative path-integral approaches, only restricted statistics are needed to obtain the control inputs. Because the Koopman operator gives a direct time-evolution of observables, it could have computational merits avoiding any sampling stages. Hence, it is valuable to seek theoretical connections of the Koopman operator with various research topics such as control issues.

The present paper explores a connection among the stochastic HJB equation, Kappen's path-integral approach, and the Koopman operator. We start from a control problem expressed in stochastic differential equations. Then, a discrete-state system is derived via the usage of the Feynman-Kac formula, It\^o's lemma, and basis expansions. The numerical solution for the derived discrete-state system gives feedback control inputs for the original stochastic differential equation. Note that there are various types of equations for the control problems. The stochastic HJB equation is solved numerically as a partial differential equation, and the path-integral and the iterative path-integral approaches need Monte-Carlo samplings. In the present paper, coupled ordinary differential equations for a discrete-state system are derived; the key of the derivation is a basis expansion focusing on the key statistic for the control problem. Although the Koopman operator gives a map on the function space, there is no need to evaluate arbitrary functions. That is, only the statistic related to the terminal cost function is necessary, and this characteristic derives the discrete-state system naturally. As a demonstration, a control problem for the van der Pol system with noise is given and solved.

The remainder of this paper is composed as follows. Section~\ref{sec_preliminaries} gives some basics of stochastic control problems. The known connection between the stochastic HJB equation and the path-integral control is also explained. Section~\ref{sec_Koopman} gives the main contribution of the present paper; a further connection with the Koopman operator is discussed, and an explicit algorithm to evaluate feedback control inputs is given. The discussion is demonstrated for a nonlinear stochastic system, which helps the understandings of the connections and derivations. Section~\ref{sec_conclusion} gives conclusions and remaining tasks.

\section{Preliminaries}
\label{sec_preliminaries}

This section gives a brief explanation for stochastic optimal control for readers unfamiliar with control issues. As for a few more details of the derivation of the stochastic HJB equation, see Appendix~\ref{sec_appendix_HJB} and Kappen's original paper \cite{Kappen2005a}. The discussions in previous works are also reviewed.

\begin{figure}
\begin{center}
\includegraphics[width=80mm]{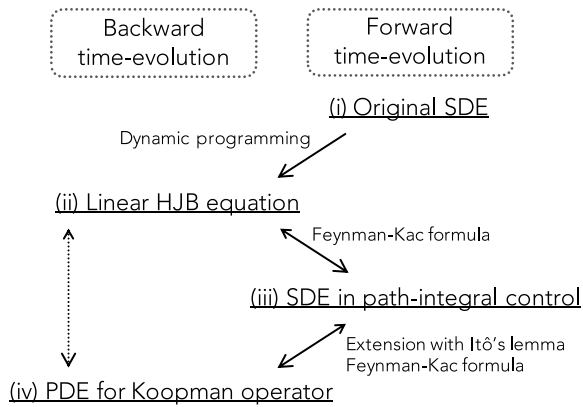}
\end{center}
\caption{An overview of the connection among the (stochastic) linear Hamilton-Jacobi-Bellman (HJB) equation, the path-integral control approach by Kappen \cite{Kappen2005a}, and the Koopman operator. Starting from the original stochastic differential equation (SDE), various types of equations are derived.
}
\label{fig_framework}
\end{figure}

Figure~\ref{fig_framework} shows an overview of the connections among the methods discussed in the present paper. Starting from the original stochastic differential equation (SDE), various types of equations are derived. Part (iv) in Fig.~\ref{fig_framework} will be explained later in Sec.~\ref{sec_Koopman}.

\subsection{Stochastic optimal control}
\label{subsec_optimal_control_intro}

Let $\bm{x}(t)$ be an $N$-dimensional stochastic variable to express the system state. For the sake of brevity of notation, sometimes the explicit time-dependency is omitted, and $\bm{x}$ is used. Consider the following stochastic differential equation for control problems:
\begin{align}
d\bm{x}  = \left(\bm{a}(\bm{x},t) + U \bm{u}(t) \right) dt + B d \bm{W}(t),
\label{eq_original_SDE}
\end{align}
where $\bm{a}(\bm{x},t)$ is an $N$-dimensional drift coefficient vector, $\bm{u}(t)$  is an $N_\mathrm{inp}$-dimensional vector for control inputs, $U$ is an $N\times N_\mathrm{inp}$ matrix related to the control setting, $B$ is an $N\times N_\mathrm{W}$ coefficient matrix for diffusion, and $d\bm{W}(t)$ stems from a conventional Wiener process with $N_\mathrm{W}$ dimensions. Note that the Wiener process satisfies the following independent property:
\begin{align}
dW_i(t) dW_j(t) = \delta_{ij} dt.
\end{align}
There is a comment for the notation of $\bm{u}(t)$. As we will see later, feedback-type control inputs are obtained, and hence $\bm{u}$ depends on $\bm{x}$. The dependency on $\bm{x}$ is abbreviated here.

Suppose a trajectory cost function $S([\bm{x}], [\bm{u}])$ consists of a terminal cost $\varphi( \bm{x}(t_\mathrm{f})): \mathbb{R}^N \to \mathbb{R}$, a running cost $V(\bm{x},t): \mathbb{R}^N \times \mathbb{R} \to \mathbb{R}$, and a quadratic control cost as follows:
\begin{align}
&S([\bm{x}], [\bm{u}]) \nonumber \\
&= \varphi( \bm{x}(t_\mathrm{f}) )
 + \int_{t_\mathrm{i}}^{t_\mathrm{f}} dt 
\left( V(\bm{x}(t),t)  + \frac{1}{2} \bm{u}(t)^\mathrm{T} R \bm{u}(t) \right),
\label{eq_trajectory_cost}
\end{align}
where $t_\mathrm{i}$ is an initial time and $t_\mathrm{f}$ is a final time. $R \in \mathbb{R}^{N_\mathrm{inp} \times N_\mathrm{inp}}$ is a weight matrix for the quadratic cost. Assume that a trajectory of $\bm{x}(t)$, $[\bm{x}]$, starts from a specific initial state $\bm{x}_\mathrm{i}$, i.e., $\bm{x}(t=t_\mathrm{i}) = \bm{x}_\mathrm{i}$.

Note that the trajectory cost function in Eq.~\eqref{eq_trajectory_cost} is defined only for a trajectory. The system is a stochastic one, and hence the following cost function should be considered:
\begin{align}
C(\bm{x}_\mathrm{i}, t_\mathrm{i}, [\bm{u}]) 
=
\mathbb{E}_{\bm{x}_\mathrm{i}} \left[ S([\bm{x}],[\bm{u}]) \right],
\end{align}
where $\mathbb{E}_{\bm{x}_\mathrm{i}} [\cdot ]$ means an expectation on trajectories with the initial condition $\bm{x}(t=t_\mathrm{i}) = \bm{x}_\mathrm{i}$.

For a time $t \in [t_\mathrm{i}, t_\mathrm{f})$, the optimal cost-to-go function is defined as
\begin{align}
J(\bm{x},t) = \min_{\bm{u}(t\to t_\mathrm{f})} C\left(\bm{x},t, \bm{u}(t \to t_\mathrm{f}) \right),
\label{eq_optimal_cost_to_go_function}
\end{align}
where $\bm{u}(t \to t_\mathrm{f})$ denotes the sequence of controls $\bm{u}(\cdot)$ on the time interval $[t, t_\mathrm{f}]$. As explained in Appendix~\ref{sec_appendix_HJB}, the optimal control inputs are evaluated via the optimal cost-to-go function as follows:
\begin{align}
\bm{u} = - R^{-1} U^\mathrm{T} \partial_{\bm{x}} J(\bm{x},t).
\label{eq_evaluate_u}
\end{align}

\subsection{Linear HJB equation}
\label{subsec_HJB}

In order to solve the stochastic optimal control problem in Sec.~\ref{subsec_optimal_control_intro}, Kappen derived a linear partial differential equation~\cite{Kappen2005a}. Kappen's derivation employs the following restriction among the coefficient matrix $B$, the matrix related to control inputs $U$, and the weight matrix for the quadratic cost:
\begin{align}
B B^\mathrm{T} = \lambda U R^{-1} U^\mathrm{T}.
\label{eq_for_lambda}
\end{align}
Note that the restriction also determines the value of $\lambda$. The restriction is necessary to obtain a linear HJB equation. As explained in Kappen's paper \cite{Kappen2005a}, no control is allowed in noiseless coordinates. Sometimes one may have noisy coordinates on which no direct controls are available. In such cases, one should choose controlled stochastic models carefully. That is, it is necessary to remove some noise terms for the control models. An example is demonstrated in Sec.~\ref{sec_Koopman}; when the noise term in the uncontrolled coordinate is small enough, the control framework works well.

When we define a function $\psi(\bm{x},t)$ via
\begin{align}
J(\bm{x},t) = - \lambda \log \psi(\bm{x},t),
\label{eq_psi_to_J}
\end{align}
the time-evolution equation for $\psi(\bm{x},t)$ is given by
\begin{align}
&\partial_t \psi(\bm{x},t) = - \mathcal{H} \psi(\bm{x},t), 
\label{eq_HJB_original}
\end{align}
where
\begin{align}
&\mathcal{H} = 
- \frac{V(\bm{x},t)}{\lambda} + \sum_{i} a_i(\bm{x},t) \partial_{x_i} 
+ \frac{1}{2} \sum_{i,j} \left[ B B^\mathrm{T}\right]_{ij} 
\partial_{x_i} \partial_{x_j}.
\label{eq_Hamiltonian_in_HJB}
\end{align}
The derivation is based on the dynamic programming method; see the original paper~\cite{Kappen2005a} and Appendix~\ref{sec_appendix_HJB}. Equation~\eqref{eq_HJB_original} is the linearized version of the stochastic HJB equation, which is the starting point of the following discussions. Note that Equation~\eqref{eq_HJB_original} is a kind of backward time-evolution one, and it is necessary to solve it in a time-reversed manner (i.e., $t_\mathrm{f} \to t_\mathrm{i}$). Additionally, note that the `initial' condition is 
\begin{align}
\psi(\bm{x},t_\mathrm{f}) = \exp(-\varphi(\bm{x})/\lambda).
\end{align}
Hence, the optimal cost-to-go function $J(\bm{x},t)$ is evaluated from $\psi(\bm{x},t)$ using Eq.~\eqref{eq_psi_to_J}. Then, we finally obtain the feedback control inputs $\bm{u}_\mathrm{opt}(\bm{x}, t)$ using Eq.~\eqref{eq_evaluate_u}.

The above discussion corresponds to the connection from (i) to (ii) in Fig.~\ref{fig_framework}. Here, two comments are given for later use.

The first comment is about the final expression. Note that the derived time-evolution equation \eqref{eq_HJB_original} is rewritten in the following form:
\begin{align}
\frac{\partial \psi}{\partial t} + \mathcal{L} \psi + g \psi = 0,
\label{eq_backward_Kolmogorov}
\end{align}
where 
\begin{align}
&\mathcal{L} = 
\sum_{i} a_i(\bm{x},t) \partial_{x_i} 
+ \frac{1}{2} \sum_{i,j} \left[ B B^\mathrm{T}\right]_{ij} 
\partial_{x_i} \partial_{x_j}, \nonumber \\
&g(\bm{x},t) = 
- \frac{V(\bm{x},t)}{\lambda}.
\end{align}
We will see later that this reformulation is beneficial for the usage of the Feynman-Kac formula.

The second comment is related to the numerical approach for the derived equation. In order to solve the partial differential equation \eqref{eq_backward_Kolmogorov} in a time-reversed manner, it is typical to employ a space-discretization method. That is, we split the continuous coordinates into a discrete lattice with a small lattice-spacing, e.g., $\Delta x \ll 1$. Then, the partial differential equation with the continuous coordinates is approximately solved via the coupled ordinary differential equations for the discrete lattice. Of course, the lattice-spacing $\Delta x$ should be small enough to obtain accurate results, and the smaller lattice-spacing leads to a larger number of equations for the coupled ordinary equations for the discrete lattice.

\subsection{SDE in path-integral control}
\label{subsec_Kappen}

The connection from (ii) to (iii) in Fig.~\ref{fig_framework} is reviewed here. In Kappen's original papers \cite{Kappen2005,Kappen2005a}, stochastic differential equations with particle-annihilation processes were derived via the Dirac's delta function and a path-integral formulation. Although the essential part is the same, here a slightly different derivation of the stochastic differential equations is given. The derivation is based on the Feynman-Kac formula, as discussed in Ref.~\cite{Thijssen2015}. As for the Feynman-Kac formula, for example, see Ref.~\cite{Kloeden_book}. As we will see here, the Feynman-Kac formula gives the stochastic differential equations in the path-integral control directly without any complicated discussions.

The Feynman-Kac formula indicates the following facts: If a partial differential equation has the form with Eq.~\eqref{eq_backward_Kolmogorov}, the solution is given as the conditional expectation
\begin{align}
\psi(\widetilde{\bm{x}},t) =
\mathbb{E}\left[
 \left. f(\bm{x}(t_\mathrm{f})) \exp \left( \int_t^{t_\mathrm{f}} g(\bm{x},\tau) d\tau \right) \right| \bm{x}(t) = \widetilde{\bm{x}}
\right],
\label{eq_Feynman_Kac}
\end{align}
where $f(\bm{x})$ is the initial condition of Eq.~\eqref{eq_backward_Kolmogorov}, i.e.,
\begin{align}
f(\bm{x}) = \exp\left( - \frac{\varphi(\bm{x})}{\lambda}  \right).
\end{align}
Note that the expectation is taken over trajectories $\bm{x}$ driven by the following stochastic differential equation
\begin{align}
d\bm{x}  =\bm{a}(\bm{x},t)  dt + B d \bm{W}(t).
\label{eq_Kappen_SDE}
\end{align}
Equation~\eqref{eq_Kappen_SDE} has the initial condition $\bm{x}(t) = \widetilde{\bm{x}}$ and develops in $[t, t_\mathrm{f}]$.

To evaluate the conditional expectation in Eq.~\eqref{eq_Feynman_Kac}, the Monte Carlo simulation is employed. Note that the particle-annihilation process is introduced in the original Kappen's papers \cite{Kappen2005,Kappen2005a}. That is, the particle is sometimes taken out of the simulation. In the formulation with the Feynman-Kac formula, the weight factor $\exp \left( \int_t^{t_\mathrm{f}} g(\bm{x},\tau) d\tau \right)$ in Eq.~\eqref{eq_Feynman_Kac} plays the role of the particle-annihilation process.

The direct Monte Carlo method for the stochastic differential equation in Eq.~\eqref{eq_Feynman_Kac} gives the control inputs, as discussed in Ref.~\cite{Kappen2005a}. Hence, there is no need to employ the space-discretization as for the HJB equation. However, the path-integral control does not give all the feedback-control inputs at once; the framework gives the control input only for a specific initial position $\bm{x}$. If one wants to evaluate the feedback control in advance as in the case of the HJB method, the Monte Carlo simulations in the path-integral control should be performed repeatedly for various initial positions. This means that the space-discretization for the initial positions is needed eventually.

\section{Koopman Operator and Control Demonstration}
\label{sec_Koopman}

This section gives the main contribution of the present paper. The connection from (iii) to (iv) in Fig.~\ref{fig_framework} is given. After general discussions, a concrete example, the stochastic van der Pol equation, is used to demonstrate the control; the explicitly derived equations will be useful to understand the general discussions.

\subsection{Derivation of PDE for Koopman operator}
\label{subsec_Koopman}

Here, a viewpoint of the Koopman operator is newly introduced. As discussed below, this viewpoint gives us the fact that we only need to focus on a specific statistic. Hence, a basis expansion is naturally introduced, which gives the coupled ordinary differential equations. After some general discussions, an explicit example for the van der Pol system will be given; the example will help the understanding of the discussion here.

Firstly, a brief introduction of the Koopman operator is given; for details, see Refs.~\cite{Williams2015, Mauroy2020}. Consider a dynamical system (with stochasticity) $(\mathcal{M}, t, F)$, where $\mathcal{M} \subseteq \mathbb{R}^N$ is the state space, $t$ is a time, and $F: \mathcal{M} \to \mathcal{M}$ is the time-evolution operator. Then, we have
\begin{align}
\bm{x}(t_\mathrm{f}) = F(\bm{x}(t_\mathrm{i}); \bm{\omega}),
\end{align}
where $\bm{\omega} \in \Omega$ is an element in the probability space associated with the dynamics. That is, the time-evolution operator $F$ changes the system state $\bm{x}$ at the initial time $t_\mathrm{i}$ to the system state $\bm{x}$ at time $t_\mathrm{f}$. Next, consider an observable $\zeta: \mathcal{M} \to \mathbb{C}$. The Koopman operator $\mathcal{K}$ associated with the time-evolution operator $F$ is defined through the composition
\begin{align}
\mathcal{K} \zeta = \zeta \circ F.
\end{align}
The Koopman operator has been used for various applications ranging from control problems, as introduced in Sec.~\ref{sec_introduction}, to prediction and change-point detection tasks \cite{Kawahara2016,Takeishi2017}.

Note that the Koopman operator gives a map for the functions of the state space. It would be typical to use some statistical quantities as $\zeta$. It indicates that we consider the time-evolution of the statistical quantities instead of that in the original state space. It is a well-known fact that the Koopman operator is \textit{linear} even when the original dynamics $F$ is nonlinear. On the other hand, the Koopman operator is infinite-dimensional. Hence, we will need finite-dimensional approximations in its practical use.

\begin{center}
\begin{table*}[tb]
\caption{Summary of the differences among the HJB equation, the path-integral control, and the proposed method.}
\label{table_summary}
\begin{tabular}{c|ccc}
\hline\hline
 & HJB & Path-integral control & Proposed \\
\hline
Equation & PDE & SDEs & Coupled ODEs\\
Solver & \quad Deterministic (e.g., Runge-Kutta) & Monte Carlo & \quad Deterministic (e.g., Runge-Kutta) \\
Discretization & Time and space & Time & Time (and basis expansions with a finite cut-off) \\
\hline\hline
\end{tabular}
\end{table*}
\end{center}

Secondly, let us continue to rewriting the path-integral control; this naturally leads to the framework from the viewpoint of the Koopman operator. Note that Eq.~\eqref{eq_Kappen_SDE} gives the time-evolution in the state space. Additionally, Eq.~\eqref{eq_Feynman_Kac} indicates that the statistic $\exp(-\varphi(x)/\lambda)$ should be focused on after the time-evolution by Eq.~\eqref{eq_Kappen_SDE}. Hence, the following definition is introduced here:
\begin{align}
z = \exp\left( - \frac{\varphi(\bm{x})}{\lambda}\right).
\end{align}
That is, the factor $f(\bm{x}) $ in Eq.~\eqref{eq_Feynman_Kac} is replaced with a new stochastic variable $z$. Then, using It\^o's lemma, we have
\begin{align}
dz = \widetilde{a}(\bm{x},z,t) dt + \widetilde{b}(\bm{x},z,t)^\mathrm{T} d\bm{W}(t),
\end{align}
where $\widetilde{a}(\bm{x},z,t)$ and $\widetilde{b}(\bm{x},z,t)^\mathrm{T}$ are adequately derived from It\^o's lemma (As for the It\^o's lennma, see also Appendix~\ref{sec_appendix_Ito_formula}.) Note that the diffusion coefficient depends on the state variable $\bm{x}$, although the constant diffusion coefficients are used in the original differential equation \eqref{eq_original_SDE} and the stochastic differential equation \eqref{eq_Kappen_SDE} in the path-integral control. Putting these stochastic variables together, we have the following \textit{extended} stochastic differential equation:
\begin{align}
d\begin{bmatrix}
\bm{x} \\
z
\end{bmatrix}
= \begin{bmatrix}
\bm{a}(\bm{x},t) \\
\widetilde{a}(\bm{x},z,t)
\end{bmatrix} dt
+ \begin{bmatrix}
B & 0 \\
\widetilde{b}(\bm{x},z,t)^\mathrm{T} & 0
\end{bmatrix}
\begin{bmatrix}
d\bm{W}(t) \\
0
\end{bmatrix}.
\label{eq_extended_SDE}
\end{align}
Let us rewrite this as
\begin{align}
d\bm{x}' = \bm{a}'(\bm{x}',t) dt + B'(\bm{x}',t) d\bm{W}'(t),
\end{align}
where $\bm{x}'$ is the $N+1$ dimensional vector including the new variable $z$, $\bm{a}'(\bm{x}',t)$ is the $N+1$ dimensional vector, $B'(\bm{x}',t)$ is the $(N+1)\times (N_\mathrm{W}+1)$ matrix, and $\bm{W}'(t)$ is the $N_\mathrm{W}+1$ vector. Then, using the Feynman-Kac formula inversely, we obtain the following backward time-evolution equation
\begin{align}
\frac{\partial \psi'}{\partial t} + \mathcal{L}' \psi' + g \psi' = 0,
\label{eq_new_backward}
\end{align}
where 
\begin{align}
\mathcal{L}'
=& \sum_{i} a'_i(\bm{x}',t) \partial_{x'_i}  \nonumber \\
&+ \frac{1}{2} \sum_{i,j} \left[ B'(\bm{x}',t) \left(B'(\bm{x}',t)\right)^\mathrm{T}\right]_{ij} \partial_{x'_i} \partial_{x'_j}.
\end{align}
Note that Eq.~\eqref{eq_new_backward} is solved backwardly $t_\mathrm{f} \to t_\mathrm{i}$, and its initial condition is
\begin{align}
\psi'(\bm{x}',t_\mathrm{f}) = \exp\left(-  \frac{\varphi(\bm{x})}{\lambda} \right).
\label{eq_new_initial_condition}
\end{align}
Because Eq.~\eqref{eq_new_backward} does not satisfy the probability conservation law, the function $\psi'(\bm{x}',t)$ is not a probability density function. Instead of the time-evolution of the system state $\bm{x}'$, we perform the time-evolution of the function $\psi'(\bm{x}',t_\mathrm{f})$ in Eq.~\eqref{eq_new_initial_condition}, i.e., 
\begin{align}
\frac{\partial \psi'}{\partial t} = \left( - \mathcal{L}' - g \right) \psi'.
\label{eq_new_backward_2}
\end{align}
The time-evolution operator $-\mathcal{L}' - g$ governs the time-evolution for the function of the state space. Hence, recalling the explanation for the Koopman operator, the integration with $-\mathcal{L}' - g$ plays a role as the Koopman operator here. Note that Eq.~\eqref{eq_new_backward_2} is linear, although the extended stochastic differential equation in Eq.~\eqref{eq_extended_SDE} is nonlinear.

\subsection{Derivation of coupled ordinary differential equations}
\label{subsec_derivation_ODEs}

We focus on the key fact that there is no need to consider a time-evolution of arbitrary functions here; the initial condition is fixed as Eq.~\eqref{eq_new_initial_condition}. This fact naturally gives the following type of basis expansion:
\begin{align}
\psi'(\bm{x}',t) =
\sum_{n_{1}, \dots, n_{N}, n_{z}} P(n_{1}, \dots, n_{N}, n_{z}, t)
x_1^{n_{1}} \dots x_N^{n_{N}} z^{n_{z}}.
\label{eq_function_expansion}
\end{align}
The basis expansion gives the initial condition for $P(n_{1}, \dots, n_{N}, n_{z}, t_\mathrm{i})$ as
\begin{align}
P(n_{1}, \dots, n_{N}, n_{z}, t_\mathrm{f}) = 
\delta_{n_1,0} \dots \delta_{n_N,0} \delta_{n_z,1}.
\end{align}
That is, it is enough to focus only on $z$ at the initial time. As for the system state $n_1, \dots, n_N$, we could employ other types of polynomials, such as the Hermite polynomials, as discussed in Ref.~\cite{Ohkubo2019}.

Then, the partial differential equation \eqref{eq_new_backward} is replaced with the coupled ordinary differential equations for $P(n_{1}, \dots, n_{N}, n_{z}, t)$. The time-evolution is also performed backwardly, i.e., $t_\mathrm{f} \to t_\mathrm{i}$. Various numerical solvers are available here; for example, the famous Runge-Kutta method is employed in the demonstration in Sec.~\ref{subsec_numerical_demonstration}. Then, the coefficients $P(n_{1}, \dots, n_{N}, n_{z}, t_\mathrm{i})$ are evaluated numerically. Hence, Eq.~\eqref{eq_evaluate_u} immediately gives the feedback control inputs at \textit{the initial time} $t_\mathrm{i}$ \textit{of the original control system}. Using the basis expansion in Eq.~\eqref{eq_function_expansion}, we have
\begin{align}
u(\bm{x},t_\mathrm{i}) = 
\frac{\lambda R^{-1} U^\mathrm{T}  \partial_{\bm{x}} \psi'(\bm{x},t_\mathrm{i})}
{
{\displaystyle\sum_{n_{1}, \dots, n_{N}, n_{z}} }
P(n_{1}, \dots, n_{N}, n_{z}, t_\mathrm{i})
x_1^{n_{1}} \dots x_N^{n_{N}} z^{n_{z}}
}.
\label{eq_u_dual}
\end{align}
The derivative of $\psi'(\bm{x},t_\mathrm{i})$ with respect to $\bm{x}$ is easily calculated because the factors $\{P(n_{1}, \dots, n_{N}, n_{z}, t_\mathrm{i})\}$ in Eq.~\eqref{eq_function_expansion} are already known. Of course, we must pay attention to the fact that $z$ is also a function of $\bm{x}$, and its derivatives should be considered adequately.

One may think that it is possible to apply the basis expansion like the one here directly to the stochastic HJB equation \eqref{eq_backward_Kolmogorov}. However, we have noticed the usage of the basis expansion through the discussions based on the Koopman operator. The expansion with respect to the desired observable $z$ is naturally obtained from the perspective of the Koopman operator, as we have seen so far.

There is a remaining comment. When the original stochastic differential equation has terms with non-polynomial functions such as $\sin$ and $\cos$, it is not straightforward to obtain the time-evolution equation with discrete-states like $P(n_{1}, \dots, n_{N}, n_{z}, t)$. We need a further variable transformation in these cases. Such a variable transformation was discussed in Ref.~\cite{Ohkubo2020}. After the variable transformation, the usage of It\^o's lemma will recover the current discussion.

A quick summary is given hare. For the stochastic optimal control, only a limited type of expectation is enough to obtain the control inputs. We can interpret the calculation of the expectation as a mapping of functions to functions, which leads to the natural introduction of the Koopman operator. The new variable $z$ is introduced via the path-integrals and It\^o's lemma, and the introduction of the new variable plays an important role in the basis expansion for $\psi'(\bm{x}',t)$. As a result, the connection revealed in the present paper gives a new numerical method to evaluate control inputs. The new numerical method is different from the previous ones. The differences are summarized in Table~\ref{table_summary}. That is, the original HJB equation is the partial differential equation (PDE); the Monte Carlo method for stochastic differential equations (SDEs) is employed in the path-integral control. In contrast, coupled ordinary differential equations (ODEs) for discrete state space are used here. Of course, as denoted in the end of Sec.~\ref{subsec_HJB}, the HJB equation also gives the coupled ordinary differential equations after the space-discretization. However, the meanings of the discrete state space is different from the proposed one; in the proposed method, the coupled ordinary differential equations for the coefficients in the basis expansion in Eq.~\eqref{eq_function_expansion} are solved numerically. The basis function is intrinsically infinite one, and hence a finite cut-off is necessary. We will see that the small cut-off is enough to reproduce the HJB results in the explicit example below.

\subsection{Demonstration for the derivation}
\label{subsec_demonstration}

\begin{figure}[b]
\begin{center}
\includegraphics[width=90mm]{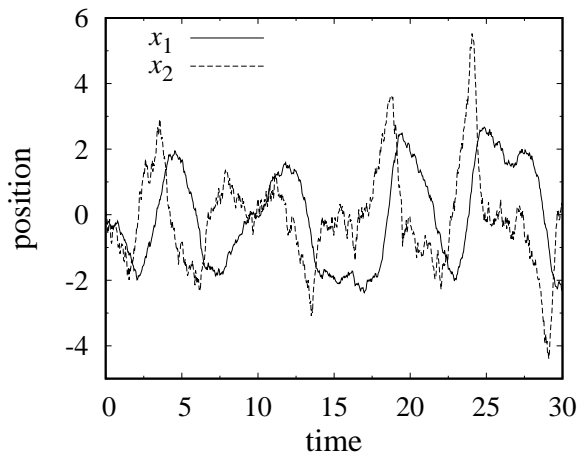}
\end{center}
\caption{Sample trajectories for Eq.~\eqref{eq_van_der_pol_original}.
There are noise terms for $x_1$ and $x_2$.
The oscillating behavior is typical for the van der Pol system.
}
\label{fig_original_model}
\end{figure}

As an example for demonstration, the following van der Pol system with noise is used here:
\begin{align}
\begin{cases}
d x_1 = x_2 dt + B_{11} dW_1(t), \\
d x_2 = \left[ \epsilon \left(1-x_1^2\right) x_2 - x_1 \right] dt + B_{22} dW_2(t).
\end{cases}
\label{eq_van_der_pol_original}
\end{align}
The van der Pol system with noise has already been used for a test of the filtering problem \cite{Lakshmivarahan2009,Ohkubo2015}.
Figure~\ref{fig_original_model} shows sample trajectories with parameters $\epsilon = 1.0$, $B_{11} = 0.1$, and $B_{22} = 1.0$. The trajectories were obtained by simulation with the Euler-Maruyama approximation \cite{Kloeden_book,Gardiner_book}; the time-discretization $\Delta t = 10^{-4}$ was employed.

\subsubsection{Dynamics and targets for control}

Here, assume that we can control only the second variable $x_2$ in Eq.~\eqref{eq_van_der_pol_original}. That is, 
\begin{align}
U = \begin{bmatrix} 0 & 0 \\ 0 & U_{22}\end{bmatrix}, \quad R = \begin{bmatrix} 0 & 0 \\ 0 & R_{22}\end{bmatrix}.
\end{align}
Hence, the following stochastic differential equation is used for our control problem:
\begin{align}
\begin{cases}
d x_1 = x_2 dt, \\
d x_2 = \left[ \epsilon \left(1-x_1^2\right) x_2 - x_1 +  U_{22} u_2(t) \right] dt + B_{22} dW_2(t),
\end{cases}
\label{eq_van_der_pol_controlled}
\end{align}
where $u_2(t) = - R_{22}^{-1} U_{22} \partial_{x_2} J(\bm{x},t)$. Note that $B_{11}$ in Eq.~\eqref{eq_van_der_pol_original} is removed; as mentioned in Sec.~\ref{subsec_optimal_control_intro}, the noise term is not permitted in the control framework here when the stochastic variables are not under direct control settings. It would be a natural assumption that $x_2$ is subject to large noise because the variable is connected to the external control inputs.

In the present work, the control targets are set as $x_1 \to 1$ and $x_2 \to 0$. Hence, the following terminal cost and running cost are used:
\begin{align}
\varphi(\bm{x}) = V(\bm{x},t) = \frac{(x_1 - x_1^\mathrm{c})^2}{2 \sigma_1^2}
+ \frac{(x_2 - x_2^\mathrm{c})^2}{2 \sigma_2^2},
\end{align}
where $x_1^\mathrm{c} = 1.0$ and $x_2^\mathrm{c} = 0.0$. In the numerical demonstration later, the other parameters are set as $\sigma_1 = 0.5$, $\sigma_2 = 0.5$, $U_{22} = 1.0$, and $R_{22} = 0.25$; these parameters lead to $\lambda = 0.25$ via Eq.~\eqref{eq_for_lambda}.

Some comments for the choice of the example are given here. There are some works for the control problem based on the data-driven Koopman operator methods. However, most of them focus on deterministic systems. Although the Ornstein-Uhlenbeck process with noise is employed as one of the examples in Ref.~\cite{Klus2020}, the system has only one coordinate. The above van der Pol system with noise has two coordinates, and the oscillating behavior makes the control problem difficult. Although the present paper aims to investigate the theoretical connection behind the control problems, it would be preferable to show the proposed framework can control the oscillating behavior. Hence, the oscillating system with noise is employed here.

\subsubsection{Coupled ordinary differential equations}

The discussion in Sec.~\ref{subsec_Koopman} gives the time-evolution operator $\mathcal{L}'$ in Eq.~\eqref{eq_new_backward}. For the example here, the following time-evolution operator is derived:


\begin{align}
&\mathcal{L}' =
x_2 \partial_{x_1}
+ (\epsilon (1-x_1^2)x_2 -x_1 )\partial_{x_2} \nonumber \\
&+ \left( 
\frac{B_{22}^2}{2} \left( - \frac{z}{\lambda \sigma_2^2} 
+ \frac{z\left(x_2 - x_2^\mathrm{c}\right)^2}{2 \lambda^2 \sigma_2^4} \right) \right. \nonumber \\
&\left. \qquad - \frac{z (x_2 - x_2^\mathrm{c})(\epsilon (1-x_1^2)x_2 -x_1 )}{\lambda \sigma_2^2}
- \frac{x_2 z (x_1 - x_1^\mathrm{c})}{\lambda \sigma_1^2}
\right) \partial_{z} \nonumber \\
& + \frac{1}{2} B_{22}^2 \partial^2_{x_2}
- \left( \frac{B_{22} z(x_2-x_2^\mathrm{c})}{\lambda \sigma_2^2} \right) \partial_{x_2} \partial_z \nonumber \\
&+ \frac{1}{4} \left( \frac{B_{22} z^2(x_2-x_2^\mathrm{c})^2}{\lambda^2 \sigma_2^4} \right) \partial_{z}^2.
\end{align}
The term without partial derivatives is obtained as follows:
\begin{align}
g = - \frac{1}{\lambda} \left( \frac{(x_1 - x_1^\mathrm{c})^2}{2 \sigma_1^2}
+ \frac{(x_2 - x_2^\mathrm{c})^2}{2 \sigma_2^2}\right).
\end{align}

The basis expansion in Eq.~\eqref{eq_function_expansion} gives coupled ordinary differential equations for $P(n_{1}, \dots, n_{N}, n_{z}, t)$. Since the derivation is tedious, a simple demonstration is shown here in a pedagogical way. Assume that we have only one variable $x$; a combination of the differential operator and an expansion similar to Eq.~\eqref{eq_function_expansion} gives
\begin{align}
\partial \sum_{n=0} P(n,t) x^n &= \sum_{n=0}^\infty n P(n,t) x^{n-1}
= \sum_{n=1}^\infty n P(n,t) x^{n-1} \nonumber \\
&= \sum_{n=0}^\infty (n+1) P(n+1,t) x^n.
\end{align}
Hence, the differential operator $\partial$ corresponds to $(n+1) P(n+1,t)$, i.e., $\partial \to (n+1) P(n+1,t)$. In similar ways, we have the following relations:
\begin{align}
\begin{cases}
&x \to P(n-1,t), \\
&x^2 \to P(n-2,t), \\
&x \partial \to n P(n,t), \\
& \partial^2 \to (n+2)(n+1)P(n+2,t), \\
&x^2 \partial \to (n-1) P(n-1,t), \\
&  x^2 \partial^2 \to n(n-1)P(n,t).
\end{cases}
\end{align}

As for the example for the demonstration, the similar discussion leads to the following coupled ordinary differential equations:
\begin{align}
\frac{\partial}{\partial t}  & P(n_1,n_2,n_z,t) \nonumber \\
=& 
- \frac{1}{2\lambda \sigma_{2}^2} P(n_1, n_2-2, n_z) 
+ \frac{x_2^\mathrm{c}}{\lambda \sigma_2^2} P(n_1, n_2-1, n_z) \nonumber \\
&- \frac{(x_2^\mathrm{c})^2}{2 \lambda \sigma_2^c} P(n_1, n_2, n_z)
- \frac{1}{2 \lambda \sigma_1^2} P(n_1-2, n_2, n_z) \nonumber \\
&+ \frac{x_1^\mathrm{c}}{\lambda \sigma_1^2} P(n_1-1, n_2, n_z) 
- \frac{(x_1^\mathrm{c})^2}{2 \lambda \sigma_1^c} P(n_1, n_2, n_z) \nonumber \\
&+ (n_1+1) P(n_1+1, n_2-1, n_z)
- \epsilon n_2 P(n_1-2, n_2, n_z) \nonumber \\
&+ \epsilon n_2 P(n_1, n_2, n_z)
- (n_2+1) P(n_1-1, n_2+1, n_z) \nonumber \\
& + \frac{\epsilon}{\lambda \sigma_2^2} n_z P(n_1-2, n_2-2, n_z) \nonumber \\
&- \frac{\epsilon x_2^\mathrm{c}}{\lambda \sigma_2^2} n_z P(n_1-2, n_2-1, n_z) \nonumber \\
&- \frac{\epsilon}{\lambda \sigma_2^2} n_z P(n_1, n_2-2, n_z)
+ \frac{\epsilon x_2^\mathrm{c}}{\lambda \sigma_2^2} n_z P(n_1, n_2-1, n_z) \nonumber \\
&- \frac{B_{22}^2}{2 \lambda \sigma_2^2} n_z P(n_1, n_2, n_z)
+ \frac{1}{\lambda \sigma_2^2} n_z P(n_1-1, n_2-1, n_z) \nonumber \\
&- \frac{x_2^\mathrm{c}}{\lambda \sigma_2^2} n_z P(n_1-1, n_2, n_z)
- \frac{1}{\lambda \sigma_1^2} n_z P(n_1-1, n_2-1, n_z) \nonumber \\
&+ \frac{x_1^\mathrm{c}}{\lambda \sigma_1^2} n_z P(n_1, n_2-1, n_z)
+ \frac{B_{22}^2}{2 \lambda^2 \sigma_2^4} n_z P(n_1, n_2-2, n_z) \nonumber \\
&- \frac{B_{22}^2 x_2^\mathrm{c}}{\lambda^2 \sigma_2^4} n_z P(n_1, n_2-1, n_z)
+ \frac{B_{22}^2 (x_2^\mathrm{c})^2}{2 \lambda^2 \sigma_2^4} n_z P(n_1, n_2, n_z) \nonumber \\
& + \frac{B_{22}^2}{2} (n_2+2)(n_2+1)P(n_1, n_2+2, n_z) \nonumber \\
&- \frac{B_{22}^2}{\lambda \sigma_2^2} n_2 n_z P(n_1, n_2, n_z) \nonumber \\
& + \frac{B_{22}^2 x_2^\mathrm{c}}{\lambda \sigma_2^2} (n_2+1) n_zP(n_1, n_2+1, n_z) \nonumber \\
&+ \frac{B_{22}^2}{2\lambda^2 \sigma_2^4} n_z (n_z-1) P(n_1, n_2-2, n_z) \nonumber \\
&- \frac{B_{22}^2 x_2^\mathrm{c}}{\lambda^2 \sigma_2^4} n_z (n_z-1) P(n_1, n_2-1, n_z) \nonumber \\
&+ \frac{B_{22}^2 (x_2^\mathrm{c})^2}{2\lambda^2 \sigma_2^4} n_z (n_z-1) P(n_1, n_2, n_z).
\label{eq_final_ODE}
\end{align}

Note that there is an infinite number of equations in the derived coupled ordinary differential equations. Hence, we need a finite cut-off. In the following numerical demonstration, the cut-offs with $n_i < 60$ for $n_1, n_2$ are employed. Note that the state for $n_z$ does not change in Eq.~\eqref{eq_final_ODE}. It means that it is enough to consider only the initial state for $n_z$.

\subsection{Numerical demonstration of control}
\label{subsec_numerical_demonstration}

Here, some numerical demonstrations for the concrete example in Sec.~\ref{subsec_demonstration} are given.
\subsubsection{Settings for numerical simulations}

Recall that the aim here is the control as $x_1 \to 1$ and $x_2 \to 0$. Other parameters are written in Sec.~\ref{subsec_demonstration}. Therefore, the numerical solution for Eq.~\eqref{eq_final_ODE} with the finite cut-offs gives the coefficients $P(n_1, n_2, n_z, t)$, which leads to the feedback control inputs directly via Eq.~\eqref{eq_u_dual}. Here, the 4th-order Runge-Kutta method with $\Delta t = 10^{-4}$ is employed. In all of the following simulations, we set $t_\mathrm{i} = 0$ and $t_\mathrm{f} = 0.1$.

\subsubsection{Numerical checks for the function $\psi$}

\begin{figure}[t]
\begin{center}
\includegraphics[width=80mm]{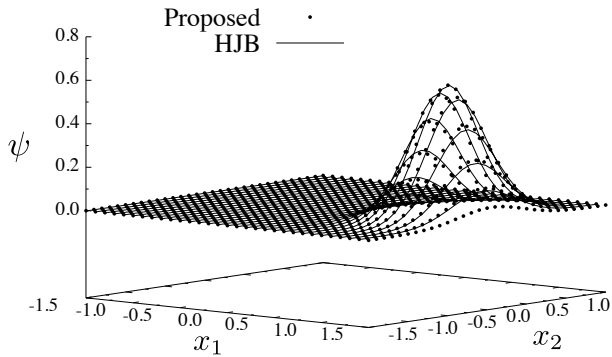}
\end{center}
\caption{A form of the evaluated $\psi$. The solid lines correspond to the numerical results for the HJB equation, and those of the proposed method are depicted with the filled circles. 
}
\label{fig_psi_results}
\end{figure}

In order to confirm that the proposed method in Secs.~\ref{subsec_Koopman} and \ref{subsec_derivation_ODEs} works well, the function $\psi$ is numerically evaluated, which is the key function to derive the feedback control inputs as discussed in Sec.~\ref{sec_preliminaries}. As explained in Sec.~\ref{subsec_derivation_ODEs}, the initial condition is set as follows:
\begin{align}
P(n_1, n_2, n_z, t = t_\mathrm{f}) =
\begin{cases}
\, 1 & \,\, \textrm{for $n_1 = 0, n_2 = 0, n_z = 1$,} \\
\, 0 & \,\, \textrm{otherwise.}
\end{cases}
\end{align}
Although the function $\psi'$ is used in Sec.~\ref{subsec_derivation_ODEs}, the substitution of an explicit value of $z$ recovers the function $\psi$. As for the comparison, the function $\psi$ was also evaluated numerically via the HJB method. In the HJB method, a region with $x_1 \in [-2, 2]$ and $x_2 \in [-2,2]$, the space-discretization with $\Delta x_1 = \Delta x_2 = 0.01$, and the time-discretization $\Delta t = 10^{-4}$ are used. Note that the HJB method needs to a little large region in order to neglect the effects from the boundary conditions. The result is shown in Fig.~\ref{fig_psi_results}, which gives a good agreement, and therefore we can see that the proposed method works well.

Note that the numbers of equations for the coupled ordinary differential equations are different intrinsically between the HJB method and the proposed one. In the HJB method, the above space-discretization needs the lattice with $(2-(-2))/0.01 \times (2-(-2))/0.01 =  400 \times 400 = 160,000$ points. In contrast, the proposed method with the above finite cut-offs gives $60 \times 60 = 3,600$ points. Note that $n_z$ does not change, and hence there is no need to consider the state-space for $n_z$. Although this is a naive comparison, the proposed method contributes to reducing the computational costs.

\subsubsection{Controlled dynamics}

Once the coefficients $P(n_1, n_2, n_z, t = t_\mathrm{i})$ are evaluated numerically, we can immediately calculate the feedback control inputs by Eq.~\eqref{eq_u_dual}. Next, the evaluated feedback control inputs are applied to the original system in Eq.~\eqref{eq_van_der_pol_original}. Note that the feedback control inputs are evaluated by removing the noise for $x_1$, as discussed in Sec.~\ref{subsec_demonstration}, but the system to be controlled has the noise term for $x_1$. We will confirm that this procedure even works well in numerical experiments.

There are some choices for the application of the obtained feedback control inputs. In the following numerical demonstration, a kind of model predictive control, i.e., the receding horizon implementation, is employed \cite{Kouvaritakis_book}. The receding horizon implementation considers a small finite time interval and calculates the feedback-control at time $0$. Note that Eq.~\eqref{eq_van_der_pol_original} is time-invariant. Hence, the evaluation for the control inputs between $t_\mathrm{i} = 0$ and $t_\mathrm{f} = 0.1$ is enough; the evaluated feedback control inputs $u(x_1,x_2) \equiv u(x_1,x_2,t=0)$ are applied all the time in the control. 

There is one additional technique from a practical viewpoint. Because of the finite cut-offs for $n_i$, large values for $x_1$ and $x_2$ may cause unexpected behavior. Hence, by defining 
\begin{align}
\widetilde{x}_i = 
\begin{cases} 
x_i^\mathrm{min} &\textrm{ if } x_i < x_i^\mathrm{min}, \\
x_i^\mathrm{max} &\textrm{ if } x_i > x_i^\mathrm{max}, \\
x_i & \textrm{otherwise},
\end{cases}
\end{align}
we employ $u(\widetilde{x}_1, \widetilde{x}_2)$ in the feedback-control. In the numerical demonstration here, the settings $x_i^\mathrm{min} = -1.5$ and $x_i^\mathrm{max} = 1.5$ for $i = 1,2$ are used.

In summary, the discussions in Sec.~\ref{subsec_Koopman} are used to evaluate the feedback control inputs, and then the feedback control inputs are applied to the original van der Pol system with noise. The control simulation is performed for the original stochastic differential equation using the Euler-Maruyama approximation with $\Delta t = 10^{-4}$. At each time step, the system state $(x_1, x_2)$ is observed, and the feedback control input $u(\widetilde{x}_1, \widetilde{x}_2)$ is evaluated. Then, we here raise the following question: Although the control input is applied only to $x_2$, is it possible to control $x_1$ together?

\begin{figure}[t]
\begin{center}
\includegraphics[width=90mm]{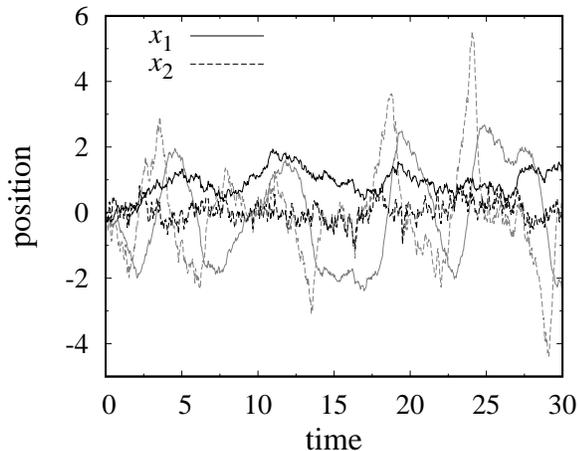}
\end{center}
\caption{Demonstration for stochastic control. The gray thin lines are original dynamics without control. The target for $x_1$ (the bold solid line) is $1$, and that for $x_2$ (the bold dashed line) is $0$.
}
\label{fig_results}
\end{figure}

Figure~\ref{fig_results} shows the numerical demonstration. The gray thin lines are the original behavior of the system, and the bold solid and dashed lines show the controlled one. The original system is noisy and oscillating, and the control input is applied only to the variable $x_2$. Despite the difficult setting, both variables are adequately controlled around the desired targets. Actually, the feedback control inputs derived from the HJB equation give a similar control result.


\section{Concluding Remarks}
\label{sec_conclusion}

The present paper clarifies one of the aspects of the importance of the Koopman operator perspectives. As pointed out in the papers for the Koopman mode decomposition \cite{Mezic2005,Williams2015}, the backward-type equation such as backward Kolmogorov equations is directly related to the evaluation for the Koopman operator. It is true even in the case of stochastic optimal controls. We can see the relationship via the path-integral approach; the target observable is considered as a new variable, and It\^o's lemma is employed. Although the Koopman operator deals with the general mapping of functions, it becomes clear that the focus only on one observable and the basis expansion are enough for our aim. Such perspective of the Koopman operators would give new insights and algorithms for other topics in physics.

Finally, some notices for practical viewpoints and further works are given. The connection revealed in the present paper leads to the coupled ordinary differential equations for the discrete state space $n_i, \dots, n_N$ and $n_z$. It will reduce the computational costs compared with the direct space-discretization for the original stochastic HJB equation, as in the example of the present paper. However, the naive application of the proposed method to high-dimensional systems would not work because of the curse of dimensionality. One may avoid the problem by using the Monte Carlo method for chemical reaction systems, i.e., the Gillespie algorithm \cite{Gillespie1977,Gillespie2007}. Additions of some terms to the derived coupled ordinary differential equations recover the probabilistic conservation law \cite{Ohkubo2013}, and hence we can employ the Monte Carlo method. Such an approach has been used for the nonlinear Kalman filter problem in Ref.~\cite{Ohkubo2015}. One may also use a discussion based on combinatorics. Recently, the methods based on combinatorics were employed to compute the Mori–Zwanzig memory integral in generalized Langevin equations \cite{Amiti2019,Zhu2020}. Similar discussions were employed to the adjoint backward equation of stochastic differential equation \cite{Ohkubo2020a}. Such methods based on combinatorics may reduce the computational cost. 
Furthermore, the applicability of the control method to stochastic chaos remains as future work. Although practical control systems usually treat non-chaotic systems, chaotic behavior in nonlinear stochastic differential equations has been studied well (e.g., see Ref.~\cite{Lin2008}), and the control problems for systems with stochastic chaos are interesting from the viewpoint of physics. In the present work, the receding horizon implementation is used, which needs only a short time evolution. Hence, the chaotic behavior would not largely affect the control because of the short time estimation. However, general discussions for control in stochastic chaos would be needed. After mitigating the curse of dimensionality problem in the Koopman approach, we will have a numerical tool to tackle this problem in the future.

The present work discusses the perspective of the Koopman operators to the stochastic optimal control. It is often enough to focus only on a specific observable quantity, and this perspective gives a new algorithm. To the best of my knowledge, this is the first attempt to employ this perspective. I hope that the discussion and the demonstration here stimulate various researchers in theoretical and practical communities.

\vspace{4mm}
\begin{acknowledgments}
This work was supported by JST, PRESTO Grant Number JPMJPR18M4, Japan.
\end{acknowledgments}

\appendix

\begin{widetext}

\section{Derivation of the Stochastic Hamilton-Jacobi-Bellman Equation}
\label{sec_appendix_HJB}

For readers' convenience, the derivation of the stochastic HJB equation is briefly explained. For details, see Kappen's original paper \cite{Kappen2005a}.

The optimal cost-to-go function $J(\bm{x},t)$ defined in Eq.~\eqref{eq_optimal_cost_to_go_function} can be rewritten as follows:
 \begin{align}
&J(\bm{x},t) \nonumber\\
&=
\min_{\bm{u}(t\to t_\mathrm{f})} 
\mathbb{E}_{\bm{x}} \left[ 
\varphi( \bm{x}(t_\mathrm{f}) )
 + 
\int_{t}^{t'} d\tau \left( \frac{1}{2} \bm{u}(\tau)^\mathrm{T} R \bm{u}(\tau) + V(\bm{x}(\tau),\tau) \right)
 + \int_{t'}^{t_\mathrm{f}} d\tau \left( \frac{1}{2} \bm{u}(\tau)^\mathrm{T} R \bm{u}(\tau) + V(\bm{x}(\tau),\tau) \right)
\right] \nonumber \\
&=
\min_{\bm{u}(t\to t')} 
\mathbb{E}_{\bm{x}} \left[ 
\int_{t}^{t'} d\tau \left( \frac{1}{2} \bm{u}(\tau)^\mathrm{T} R \bm{u}(\tau) + V(\bm{x}(\tau),\tau) \right) 
 + \min_{\bm{u}(t'\to t_\mathrm{f})} 
\mathbb{E}_{\bm{x}(t')} \left[ 
\varphi( \bm{x}(t_\mathrm{f}) )
+ \int_{t'}^{t_\mathrm{f}} d\tau \left( \frac{1}{2} \bm{u}(\tau)^\mathrm{T} R \bm{u}(\tau) + V(\bm{x}(\tau),\tau) \right)
\right]
\right] \nonumber \\
&=
\min_{\bm{u}(t\to t')} 
\mathbb{E}_{\bm{x}} \left[ 
\int_{t}^{t'} d\tau \left( \frac{1}{2} \bm{u}(\tau)^\mathrm{T} R \bm{u}(\tau) + V(\bm{x}(\tau),\tau) \right)
+ J(\bm{x}(t'), t')
\right].
\end{align}
Here, we set $t' = t + dt$ and perform the Taylor expansion for $J(\bm{x}(t'),t')$ around $t$. Taking the expansion up to the first order in $dt$, we have
\begin{align}
\mathbb{E}_{\bm{x}} \left[ J(\bm{x}(t+dt),t+dt) \right] 
&= 
\mathbb{E}_{\bm{x}} \left[ 
J(\bm{x},t) + \partial_t J(\bm{x},t) dt + \left( \partial_{\bm{x}}J(\bm{x},t) \right)^\mathrm{T} d\bm{x}
+ \frac{1}{2} \mathrm{Tr}\left( \left(\partial^2_{\bm{x}} J(\bm{x},t) \right) \left(d\bm{x} d\bm{x}^\mathrm{T}\right)\right)
\right] \nonumber \\
&= 
J(\bm{x},t) + \partial_t J(\bm{x},t) dt 
+ \left( \partial_{\bm{x}}J(\bm{x},t) \right)^\mathrm{T} \left( \bm{a}(\bm{x},t)+U\bm{u}(t) \right) dt 
+ \frac{1}{2} \mathrm{Tr} \left(  B B^\mathrm{T} \partial^2_{\bm{x}} J(\bm{x},t)  \right) dt,
\end{align}
where
\begin{align}
\mathrm{Tr}\left(\left(\partial^2_{\bm{x}} J(\bm{x}(t),t)\right) \left(d\bm{x} d\bm{x}^\mathrm{T}\right)\right) 
= \sum_{ij} \frac{\partial^2 J(\bm{x}(t),t)}{\partial x_i \partial x_j} dx_i dx_j
\end{align}
and 
\begin{align}
\mathrm{Tr}\left(  B B^\mathrm{T}  \partial^2_{\bm{x}} J(\bm{x}(t),t)\right) = \sum_{ij} \left[ B B^\mathrm{T} \right]_{ij} 
\frac{\partial^2 J(\bm{x}(t),t)}{\partial x_i \partial x_j}.
\end{align}
Hence, we have the following time-evolution equation for $J(\bm{x},t)$ with taking the limit $dt \to 0$:
\begin{align}
-\partial_t J(\bm{x},t) 
= \min_{\bm{u}} 
\left(
\frac{1}{2} \bm{u}(t)^\mathrm{T} R \bm{u}(t) + V(\bm{x},t) 
+ \left( \bm{a}(\bm{x},t)+U\bm{u}(t) \right)^\mathrm{T} \partial_{\bm{x}}J(\bm{x},t) 
 + \frac{1}{2} 
\mathrm{Tr}\left(  B B^\mathrm{T} \partial^2_{\bm{x}} J(\bm{x},t)\right) \right).
\end{align}
The minimization with respect to $\bm{u}$ gives 
\begin{align}
\bm{u} = - R^{-1} U^\mathrm{T} \partial_{\bm{x}} J(\bm{x},t),
\end{align}
which leads to 
\begin{align}
-\partial_t J(\bm{x},t) 
=& 
- \frac{1}{2} (\partial_{\bm{x}}J)^\mathrm{T} U R^{-1}U^\mathrm{T} \partial_{\bm{x}}J(\bm{x},t)
+ V(\bm{x},t)  
+ \bm{a}(\bm{x},t)^\mathrm{T}  \partial_{\bm{x}} J(\bm{x},t)
+ \frac{1}{2} 
\mathrm{Tr}\left(  B B^\mathrm{T} \partial^2_{\bm{x}} J(\bm{x},t)\right).
\label{eq_appendix_time_evolution_for_J_final}
\end{align}
The time-evolution equation \eqref{eq_appendix_time_evolution_for_J_final} should be solved with the boundary condition $J(\bm{x},t_\mathrm{f}) = \varphi( \bm{x} )$.

Next, the nonlinearity in Eq.~\eqref{eq_appendix_time_evolution_for_J_final} is removed. By defining 
\begin{align}
J(\bm{x},t) = - \lambda \log \psi(\bm{x},t),
\end{align}
where a constant $\lambda$ to be defined, we rewrite the first and final terms in r.h.s. of Eq.~\eqref{eq_appendix_time_evolution_for_J_final} as
\begin{align}
&- \frac{1}{2} (\partial_{\bm{x}}J)^\mathrm{T} U R^{-1}U^\mathrm{T} \partial_{\bm{x}}J(\bm{x},t)
+ \frac{1}{2} \mathrm{Tr}
\left( B B^\mathrm{T} \partial^2_{\bm{x}} J(\bm{x}(t),t) \right) \nonumber \\
&\quad = - \frac{\lambda^2}{2 \psi^2} \sum_{ij} \frac{\partial \psi}{\partial x_i}
\left[ U R^{-1} U^\mathrm{T} \right]_{ij} \frac{\partial \psi}{\partial x_j}
+ \frac{\lambda}{2 \psi^2} \sum_{ij} \left[ B B^\mathrm{T} \right]_{ij} 
 \frac{\partial \psi}{\partial x_i} \frac{\partial \psi}{\partial x_j}
 - \frac{\lambda}{2 \psi} \sum_{ij}  \left[ B B^\mathrm{T} \right]_{ij} 
 \frac{\partial^2 \psi}{\partial x_i \partial x_j}.
\end{align}
Hence, it is easy to see that the terms quadratic in $\psi$ vanish when we set 
\begin{align}
B B^\mathrm{T} = \lambda U R^{-1} U^\mathrm{T}.
\label{eq_appendix_for_lambda}
\end{align}
Note that Eq.~\eqref{eq_appendix_for_lambda} also determines the actual value of $\lambda$.
Finally, the following equation is derived
\begin{align}
&\partial_t \psi(x,y,t) = - \mathcal{H} \psi(x,y,t),  \\
&\mathcal{H} = 
- \frac{V(\bm{x},t)}{\lambda} + \sum_{i} a_i \partial_{x_i} 
+ \frac{1}{2} \sum_{i,j} \left[ B B^\mathrm{T}\right]_{ij} \partial_{x_i} \partial_{x_j},
\end{align}
which corresponds to Eq.~\eqref{eq_HJB_original} and Eq.~\eqref{eq_Hamiltonian_in_HJB}.

\section{Brief Summary of It\^o's Lemma}
\label{sec_appendix_Ito_formula}

As for the details of It\^o's lemma, for example, see Ref.~\cite{Gardiner_book}. Consider the following multivariate stochastic differential equation:
\begin{align}
d\bm{x} = \bm{a}(\bm{x},t) dt + B(\bm{x},t) dW(t),
\end{align}
where $\bm{a}(\bm{x},t)$ is a vector for drift terms, $B(\bm{x},t)$ is a matrix for diffusion terms. $W(t)$ is a vector for Wiener processes, whose elements satisfy
\begin{align}
dW_i(t) dW_j(t) = \delta_{ij} dt,
\end{align}
that is, $dW_i(t)$ and $dW_j(t)$ are independent each other. It\^o's lemma plays an important role to consider an arbitrary function of $f(\bm{x})$. In the stochastic calculus, a change of variables is not given by conventional differential calculus, and It\^o's lemma suggests the following formula:
\begin{align}
d f(\bm{x}) =&
\left\{
\sum_i a_i(\bm{x},t) \partial_{x_i} f(\bm{x})
+ \frac{1}{2} \sum_{i,j} \left[ B(\bm{x},t) B^\mathrm{T}(\bm{x},t) \right]_{ij} \partial_{x_i} \partial_{x_j} f(\bm{x})
\right\} dt 
+ \sum_{i,j} B_{ij}(\bm{x},t) \partial_{x_i} f(\bm{x}) \, d W_j(t).
\end{align}
The second term in the curly bracket is different from the conventional calculus.

\end{widetext}



\end{document}